\documentclass[leqno,12pt]{elsarticle}
\usepackage{amssymb}
\usepackage{amsfonts,amscd}
\usepackage{amsmath}
\usepackage{hyperref}
 \numberwithin{equation}{section}
\newtheorem{theo}{Theorem}[section]
\newtheorem{lemma}[theo]{Lemma}

\newtheorem{prop}[theo]{Proposition}
\newtheorem{remark}[theo]{Remark}

\newtheorem{example}[theo]{Example}

\newtheorem{defi}[theo]{Definition}
\newproof{pf}{Proof}
\newcommand{\inter}{\operatorname{int}}
\newcommand{\clos}{\operatorname{cl}}

\begin{document}
\title{Compact bilinear operators on asymmetric normed spaces}
\author{S. Cobza\c{s} }
\address{\it Babe\c s-Bolyai University, Department  of Mathematics,
 400 084 Cluj-Napoca, Romania,\\ E-mail: scobzas@math.ubbcluj.ro}\;\;
\begin{abstract}
The paper is concerned with compact bilinear operators on asymmetric normed spaces. The study of multilinear operators on asymmetric normed spaces was initiated by Latreche and Dahia, Colloq. Math. (2020). We go further in this direction and prove a Schauder type theorem on the compactness of the adjoint of a compact bilinear operator and   study the ideal properties of spaces of compact bilinear operators. These extend some results of Ramanujan and Schock, Linear and  Multilinear Algebra  (1985), and Ruch, ibid. (1989), on compact bilinear operators on Banach spaces.  On the space of bilinear forms one introduces the analog of the weak$^*$-topology, called the $w^2$-topology, and one proves an Alaoglu-Bourbaki type theorem -- the $w^2$-compactness of the closed unit ball.


\end{abstract}

\maketitle

\date{\today}
\maketitle

\section{Introduction}
The theory of multilinear operators on Banach spaces emerged as a counterpart to that of  linear operators and concerned mainly
extension properties, ideal theory, interpolation  and  absolutely summing multilinear operators, see, for instance, the survey papers  \cite{sanchez00}, \cite{pelleg16}, \cite{pelleg11}.

Ramanujan and Schock \cite{schock85} defined the adjoint (or conjugate) of a bilinear operator, proved a Schauder type theorem on the compactness of the adjoint of a compact bilinear  operator and studied the ideal properties of some classes of bilinear operators. Their results were completed by Ruch \cite{ruch89} (see also \cite{ruch96}).

The first attempt to study multilinear  operators on asymmetric normed spaces was made by Latreche and Dahia \cite{latreche20}. The aim of this paper is to go further in this direction and to extend the results of Ramanujan, Schock and Ruch on compact bilinear operators on Banach spaces to the asymmetric case. We define  continuous bilinear operators on asymmetric normed spaces and their adjoints, prove a Schauder type theorem on the precompactness of the adjoint of a precompact bilinear operator and study some ideal properties of spaces of precompact operators. We also introduce the analog of the weak$^*$-topology on the space of bilinear  forms (denoted by  $w^2$) and prove an Alaoglu-Bourbaki type theorem -- the $w^2$-compactness of the closed unit ball of this space.
\section{Preliminaries}
For reader's convenience, we present in this section   some  results on quasi-metric spaces, asymmetric normed spaces and their duals, and on linear operators acting between these spaces. We restrict the presentation to  notions and results used in the development of the main topic. Full details and further results can be found in the book \cite{Cobzas}.
\subsection{Quasi-metric spaces}
A {\it quasi-pseudometric} on a set $X$ is a mapping $d: X\times X\to [0,\infty)$ satisfying the following conditions:
\begin{align*}%
\mbox{(QM1)}&\qquad d(x,y)\geq 0, \quad and  \quad d(x,x)=0;\\
\mbox{(QM2)}&\qquad d(x,z)\leq d(x,y)+d(y,z),  %
\end{align*}%
for all $x,y,z\in X.$ If further
$$%
\mbox{(QM3)}\qquad d(x,y)=d(y,x)=0\Rightarrow x=y,
$$%
for all $x,y\in X,$ then $d$ is called a {\it quasi-metric}. The pair $(X,d)$ is called a {\it
quasi-pseudometric space}, respectively  a {\it quasi-metric space}\footnote{In \cite{Cobzas} the term ``quasi-semimetric" is used instead of ``quasi-pseudometric"}. The adjoint of the quasi-pseudometric
$d$ is the quasi-pseudometric $\bar d(x,y)=d(y,x),\, x,y\in X.$ The mapping $
d^s(x,y)=\max\{d(x,y),\bar d(x,y)\},\,$ $ x,y\in X,$ is a pseudometric on $X$ which is a metric if and
only if $d$ is a quasi-metric.

If $(X,d)$ is a quasi-pseudometric space, then for $x\in X$ and $r>0$ we define the balls in $X$ by the formulae %
\begin{equation}\label{def.balls-qm}\begin{aligned}
B_d(x,r)=&\{y\in X : d(x,y)<r\} \; \mbox{-\; the open ball, and }\\ %
B_d[x,r]=&\{y\in X : d(x,y)\leq r\} \; \mbox{-\; the closed ball. } %
\end{aligned} \end{equation}

 The topology $\tau_d$ (or $\tau(d)$) of a quasi-pseudometric space $(X,d)$ can be defined starting from the family
$\mathcal{V}_d(x)$ of neighborhoods  of an arbitrary  point $x\in X$:%
\begin{equation*}
\begin{aligned}
V\in \mathcal{V}_d(x)\;&\iff \; \exists r>0\;\mbox{such that}\; B_d(x,r)\subseteq V\\
                             &\iff \; \exists r'>0\;\mbox{such that}\; B_d[x,r']\subseteq V. %
\end{aligned} %
\end{equation*}

The convergence of a sequence $(x_n)$ to $x$ with respect to $\tau_d,$ called $d$-convergence and
denoted by
$x_n\xrightarrow{d}x,$ can be characterized in the following way %
\begin{equation}\label{char-rho-conv1} %
         x_n\xrightarrow{d}x\;\iff\; d(x,x_n)\to 0. %
\end{equation} %

Also
\begin{equation}\label{char-rho-conv2} %
         x_n\xrightarrow{\bar d}x\;\iff\;\bar d(x,x_n)\to 0\; \iff\; d(x_n,x)\to 0. %
\end{equation} %

As a space equipped with two topologies, $\tau_d$ and   $\tau_{\bar d}\,$, a quasi-pseudometric space can be viewed as a bitopological space in the sense of Kelly \cite{kelly63}.

The following   topological properties are true for  quasi-pseudometric spaces.
    \begin{prop}[see \cite{Cobzas}]\label{p.top-qsm1}
   If $(X,d)$ is a quasi-pseudometric space, then the followings hold.
   \begin{enumerate}
   \item[\rm 1.] The ball $B_d(x,r)$ is $\tau_d$-open and  the ball $B_d[x,r]$ is
       $\tau_{\bar{d}}$-closed. The ball    $B_d[x,r]$ need not be $\tau_d$-closed.
     \item[\rm 2.]
   The topology $\tau_d$ is $T_0$ if and only if  $d $ is a quasi-metric.   \\ The topology $\tau_d$ is $T_1$ if and only if
   $d(x,y)>0$ for all  $x\neq y$  in $X$.
      \item [\rm 3.]  For every fixed $x\in X,$ the mapping $d(x,\cdot):X\to (\mathbb{R},|\cdot|)$ is
   $\tau_d$-usc and $\tau_{\bar d}$-lsc. \\
   For every fixed $y\in X,$ the mapping $d(\cdot,y):X\to (\mathbb{R},|\cdot|)$ is $\tau_d$-lsc and
   $\tau_{\bar d}$-usc.
\end{enumerate}     \end{prop}

\subsection{Extended quasi-metric spaces}
An \emph{extended pseudometric }(\emph{metric}) (called also generalized metric) on a set $X$ is  a mapping $d:X\times Y\to[0,\infty]$ satisfying    the usual axioms of a pseudometric (metric).
Extended metric spaces $(X,d)$  were introduced by W. A. J. Luxemburg in \cite{lux1}, \cite{lux2} in connection with the method of successive approximation and fixed points. These results   were completed by A. F. Monna \cite{monna61} and M. Edelstein \cite{edelst64}.  Further results were obtained by  J. B. Diaz and B. Margolis \cite{diaz-margol68,margolis68} and C. F. K.  Jung \cite{jung69}.
For some recent results on generalized metric spaces see \cite{beer13} and \cite{czerw-krol15}, \cite{czerw-krol16a}, \cite{czerw-krol16b}, \cite{czerw-krol17}.

Recently, G. Beer and J. Vanderwerf \cite{beer15}, \cite{beer-vdw15a}  considered  vector spaces equipped with norms that  can take infinite values, called  by them  ``extended norms" (see also \cite{czerw-krol16c}).

As it was shown by Jung \cite{jung69} the relation
\begin{equation}\label{eq1.ex-qm} x\sim y\iff d(x,y)<\infty,\;\; x,y\in X,\end{equation}
is an equivalence relation on an extended pseudometric  space $(X,d)$,   decomposing it into pairwise disjoint equivalence classes $X_i,\, i\in I$. The restriction $d_i=d|_{X_i\times X_i}$ is a pseudometric on $X_i$ and many properties of $X$ can be reduced to those of the  pseudometric spaces $(X_i,d_i),\, i\in I$. For instance, an extended pseudometric  space is complete if and only if    $(X_i,d_i)$ is a complete pseudometric  space for every $i\in I$.

By analogy, an \emph{extended quasi-pseudometric}   on a set $X$ is a mapping $d:X\times X\to[0,\infty]$ satisfying the conditions (QM1), (QM2) from above. If further, (QM3) also holds, then $d$ is called an   \emph{extended quasi-metric}. In this case, due to the lack of symmetry,  the relation \eqref{eq1.ex-qm} is  only reflexive and transitive. The topology generated by an extended quasi-pseudometric   is defined as in the case of quasi-pseudometric spaces, via the balls \eqref{def.balls-qm}.

\subsection{Precompactness  in quasi-pseudometric spaces}
Recall that a subset $Y$ of a pseudometric space $(X,d)$ is called \emph{precompact} (or \emph{totally bounded}) if for every $\varepsilon>0$ there exists a finite subset $Z$ of $Y$ such that $Y\subseteq \bigcup\{B_d(z,\varepsilon): z\in Z\}$. One obtains the same notion if we take $Z\subseteq X$ and/or closed balls. Also every subset of a precompact set is  precompact.

If $(X,d)$ is only a quasi-pseudometric, then these notions are different.
\begin{defi}\label{def.prec}
  A subset $Y$ of a quasi-pseudometric  space $(X,d)$ is called:
  \begin{itemize}\item
  $d$-\emph{precompact} if for every $\varepsilon>0$ there exists a finite subset $Z$ of $Y$ such that
  \begin{equation}\label{eq1.prec}
Y\subseteq \bigcup\{B_d(z,\varepsilon): z\in Z\};\end{equation}
\item \emph{outside  $d$-precompact} if for every $\varepsilon>0$ there exists a finite subset $Z$ of $X$ such that
\eqref{eq1.prec} holds;
\item \emph{hereditarily  $d$-precompact} if   every subset of $Y$ is  $d$-precompact.
\end{itemize}\end{defi}

Similar notions can be defined for $\bar d$.
Some remarks are in order.
\begin{remark}\label{re.prec}  Let $(X,d)$ be a quasi-pseudometric space.
\begin{enumerate}\item[\rm 1.]  A $d$-precompact subset of $X$ is outside $d$-precompact. There exist outside $d$-precompact sets  which are not $d$-precompact.
\item[\rm 2.]  A $d^s$-precompact set is $d$-precompact and $\bar d$-precompact. There exist sets which are both $d$- and $\bar d$-precompact   but not $d^s$-precompact.
\item[\rm 3.] There exist $d$-precompact sets which are not hereditarily $d$-precompact. Outside $d$   ($\bar d$)-precompactness   and $d^s$-precompactness are hereditary pro\-perties.
     \end{enumerate}\end{remark}

For a discussion, relations with compactness  and proofs, see Subsection 1.2.2 of \cite{Cobzas}.  We mention only the following result proved in \cite{kunzi-reilly93} (see also \cite[Prop. 1.2.36]{Cobzas}), the analog of a well-known result in metric spaces.

\begin{prop}\label{p1.h-prec} A subset $Y$ of a quasi-pseudometric space is hereditarily precompact if an only if every sequence in $Y$ has a left $K$-Cauchy subsequence.
\end{prop}

Recall that a sequence $(x_n)_{n\in\mathbb{N}}$ in a quasi-pseudometric space is called \emph{left} (\emph{right}) $K$-\emph{Cauchy}  if for every $\varepsilon>0$ there exists $n_\varepsilon$ such that for all $m,n\in\mathbb{N}$ with $n_\varepsilon\le n\le m,\,$
   $d(x_n,x_m)<\varepsilon$ (resp.  $d(x_m,x_n)<\varepsilon$).
   \begin{remark}
Notice that the  precompactness notions given in  Definition \ref{def.prec}    can be  transposed verbatim  to  an extended quasi-metric space $(X,d)$.
   \end{remark}

\subsection{Asymmetric normed spaces}

An \emph{asymmetric norm}   on a real vector space $X$ is a functional
$p:X\to [0,\infty)$ satisfying the conditions
  \begin{equation}\label{def.asyn}\begin{aligned}
&\mbox{(AN1)}\;\; p(x)=p(-x)=0\Rightarrow x=0;\quad\mbox{(AN2)}\;\;
p(\alpha x)=\alpha p(x);\\
&\mbox{(AN3)}\;\; p(x+y)\leq p(x)+p(y)\;,
\end{aligned}\end{equation}
for all $x,y\in X$ and $\alpha \geq 0.$

In some instances, the value $+\infty$ will be allowed for $p$
in which case we shall call $p$ an \emph{extended asymmetric norm } (or \emph{seminorm}).

An asymmetric seminorm $p$ defines a
quasi-pseudometric $d_p$ on $X$ through the formula
\begin{equation}\label{def.qm}
d_p(x,y)=p(y-x),\; x,y\in X\;.
\end{equation}

If $p$ is an asymmetric norm, then $d_p$ is a quasi-metric.

 Define the conjugate asymmetric seminorm $\bar p$ and the seminorm  $p^s$ by
 \begin{equation}\label{def.bar-p}
 \bar p(x)=p(-x)\quad \mbox{and}\quad p^s(x) =\max\{p(x),p(-x)\}\;,
 \end{equation}
 for $x\in X$.

 Obviously, $p^s$ is a norm when $p$ is an asymmetric norm and $(X,p^s)$ is a normed space.
 Also,
\begin{equation*}
\bar d_p=d_{\bar p} \quad\mbox{and}\quad   d^s_p=d_{p^s}\;.
\end{equation*}

The closed unit ball of an asymmetric normed space $(X,p)$ will be denoted by $B_p$, that is
\begin{equation*}
B_p=\{x\in X:p(x)\le 1\}\,.\end{equation*}

It follows that and arbitrary closed ball satisfies the equality \begin{equation}\label{eq1.ball-qm}B_p[x,r]=x+rB_p\,.\end{equation}
 
 A subset $Y$ of an asymmetric normed space is called $p$-\emph{bounded} if $\,\sup p(Y)  <\infty.$
 
\begin{remark}\label{re.prec-bd}
   A $p$-precompact  subset  of an asymmetric normed space is  $p$-bounded.
\end{remark}

 Indeed, for $\varepsilon=1$ there exist $y_1,\dots,y_n$ in $Y$ s.t. for every $y\in Y$ there is $k_0\in \{1,\dots,n\}$  s.t. $p(y-y_{k_0})\le 1$. But then
 $$
 p(y)\le p(y-y_{k_0})+p(y_{k_0})\le 1+\max\{ p(y_k): 1\le k\le n\}\,,$$
 for all $y\in Y$, which shows that $Y$ is $p$-bounded.

The following example is essential in the study of the duality for asymmetric normed spaces.
\begin{example}\label{ex1.Ru} On $\mathbb{R}$ consider the functional $u:\mathbb{R}\to \mathbb{R}_+$ given by $u(\alpha)=\alpha^+=\max\{\alpha,0\},\; \alpha\in \mathbb{R}.$
 Then $u$ is an asymmetric norm on $ \mathbb{R}$ with conjugate  $\bar u(\alpha)=\alpha^-=\max\{-\alpha,0\}$ and the associated (symmetric) norm $u^s(\alpha)=|\alpha|$.
\end{example}

\subsection{Normed cones}

A \emph{cone} (called sometimes a \emph{convex cone}) is a subset $Z$ of a linear space $X$ such that $w+z\in Z$ and
$\lambda z\in Z$ for all $w,z\in Z$ and $\lambda \geq 0.$ A cone is also called  a \emph{semilinear space}.

\begin{remark}\label{rem.cones}
In fact, one can define an abstract notion of cone as a set $X$ with two operations, addition +  which is supposed to be commutative, associative and having a neutral element denoted by 0, and
  multiplication by nonnegative scalars (denoted by $\,\cdot$),  satisfying the properties
 \begin{equation}\label{def.cone}
 \begin{aligned}
 &\mbox{(i)}\;\; (\lambda\mu)a=\lambda(\mu a),\quad \mbox{(ii)}\;\; \lambda (a+b)=\lambda a+\lambda b, \quad
 \mbox{(iii)}\;\;
 (\lambda+\mu)a =\lambda a+\mu a,\\ &\mbox{(iv)} \;\; 1\cdot a=1,\quad\mbox{ and} \quad \mbox{(v)}\; \;0\cdot
 a=0\;.
\end{aligned}
  \end{equation}

  The cone  $X$ is called \emph{cancellative} if  $$a+c=b+c\;\Rightarrow\; a=b\,,$$ for all $a,b,c\in X.$ A cone  $X$
  is  cancellative if and only if it can be embedded in a vector space.  This means that there exist a linear space $\widetilde X$ and an injective linear mapping
  $j:X\to \widetilde X$ such that $\widetilde X=j(X)-j(X)$, i.e. the cone $j(X)$ is generating in $\widetilde X$ (see \cite[pp. 114-116]{Fuc-Luk}). Also every linear operator $A$ from $X$ to a vector space $Y$ admits a unique linear extension $\widetilde A:\widetilde X\to Y$.

  By linearity in the case of a cone we understand additivity and positive homogeneity.
      The theory of locally convex cones, with applications to Korovkin type approximation theory
  for    positive operators and to vector-measure theory, is developed in the books by Keimel and Roth
  \cite{K-Roth92} and   Roth \cite{Roth09}, respectively.
\end{remark}

We shall always suppose that the considered cones are contained in some vector spaces.

Let $X$ be a cone. A \emph{norm} on $X$ is a   positively homogeneous, subadditive mapping $p:X\to \mathbb{R}_+$   such that
\begin{equation}\label{def1.norm-cone}
x,-x\in X\;\mbox{ and }\;p(x)=0=p(-x)\;\Rightarrow\; x=0\,.\end{equation}

If the stronger condition
\begin{equation}\label{def2.norm-cone}
p(x)=0\;\Rightarrow\; x=0\,\end{equation}
 holds for all $x\in X$, then we say that $p$ is a $T_1$-\emph{norm}.

As   in the case of asymmetric normed spaces the unit ball of a normed cone $(X,p)$ will be denoted by $B_p$, that is,
\begin{equation*}
B_p=\{x\in X: p(x)\le 1\}\,.\end{equation*}

The equality \eqref{eq1.ball-qm} holds in this case too.

In \cite{romag04} an extended quasi-pseudometric $d_p$ was defined on a normed cone $(X,p)$ by
\begin{equation}\label{def.ext-qm}
d_p(x,y)=\begin{cases}
 \inf\{p(z):   z\in X,\, y=x+z\}\quad &\mbox{if }\; y\in x+X;\\
\infty\quad& \mbox{otherwise}\,.
\end{cases} \end{equation}

In our case (of a cancellative cone) this definition becomes
\begin{equation}\label{def.ext-dist}
d_p(x,y)=\begin{cases}p(z)\quad &\mbox{if }\; y=x+z\;\mbox{for some}\;z\in X;\\
\infty\quad& \mbox{otherwise}\,.
\end{cases}\end{equation}

If the cone is contained in a vector space, then
$$
y\in x+X\iff y-x\in X.$$

The corresponding topology will be denoted by $\tau(p)$ (or by $\tau_p$).
\begin{prop}[\cite{romag04}] Let $(X,p)$ be a normed cancellative cone and let $d_p:X\times X\to[0,\infty]$ be  defined by \eqref{def.ext-dist}.
Then:
\begin{itemize}\item[\rm (i)]
  the functional  $d_p$ is an  extended quasi-pseudometric, which is an extended quasi-metric if $p$ is
  a $T_1$-norm;
 \item[\rm (i)]
  the  extended quasi-pseudometric $d_p$ is translation invariant and positively homogeneous, that is,
$$
d_p(x+z,y+z)= d_p(x,y)\quad\mbox{and}\quad d_p(tx,ty)=td_p(x,y)\,$$
for all $x,y,z\in X$ and $t>0$;
\item[\rm (iii)] the balls $B_{d_p}(x,\varepsilon)$ are $\tau(p)$-open as well as the translation mapping $x\mapsto x+z,\, x\in X$, for any  $z\in X$.
\end{itemize}\end{prop}

\textbf{Note.}\, All topological and metric (e.g. precompactness) notions in a normed cone will be considered with respect to this extended metric.

\subsection{Continuous linear operators on normed cones}\label{SS.lin-ops-asyn}
A mapping $A$ between two cones $X,Y$ is called \emph{linear} if it is additive and positively homogeneous.
If $X,Y$ are vector spaces, then linearity means additivity and  homogeneity.
\begin{prop}\label{p1.lin-op-cone}
  Let $(X,p)$ and $(Y,q)$ be normed cones and let $A:X\to Y$  be  a linear mapping. Then $A$ is $\tau(p)$ to $\tau(q)$ continuous if and only if
  there exists $\beta\ge 0$ such that
\begin{equation}\label{eq1.lin-op-cone}
  q(Ax)\le\beta p(x)\,,\end{equation}
  for all $x\in X$.
\end{prop}

A number $\beta\ge 0$ for which \eqref{eq1.lin-op-cone} holds is called a \emph{semi-Lipschitz constant} for $A$.

For  a continuous linear operator $A:X\to Y$ between two normed cones $(X,p),(Y,q)$ put
\begin{equation}\label{def.norm-op-cone}
\|A|_{p,q}=\sup\{q(Ax): x\in X,\, p(x)\le 1\}\,\end{equation}
and denote by  $\mathcal L_{p,q}(X,Y)$ the set of all  continuous linear operators from $(X,p)$ to $(Y,q)$.
\begin{prop}\label{p2.lin-op-cone} Let $(X,p),(Y,q)$ be normed cones.
\begin{enumerate}\item[\rm 1.]The set $\mathcal L_{p,q}(X,Y)$ is a cone and the functional $\|\cdot|_{p,q}$
given by \eqref{def.norm-op-cone} is a norm on $\mathcal L_{p,q}(X,Y)$.
\item[\rm 2.]
The norm $\|A|_{p,q}$ of $A\in \mathcal L_{p,q}(X,Y)$ is the smallest semi-Lipschitz constant $\beta$ for which \eqref{eq1.lin-op-cone} holds.
\end{enumerate}\end{prop}

 Denote by $L(X,Y)$ the space of all linear operators between two linear spaces $X$, $Y$.
If   $(X,p)$ and $(Y,q)$ are asymmetric seminormed spaces, then  we can  consider  several spaces of linear operators.
For $\mu\in \{p,\bar p,p^s\}$ and $\nu\in \{q,\bar q,q^s\},$ the $(\mu,\nu)$-continuity of linear operators and the space
$\mathcal L_{\mu,\nu}(X,Y)$ are defined accordingly.  The space of all continuous linear operators between  the
associated  seminormed spaces $X_s=(X,p^s)$ and $Y_s=(Y,q^s)$ is denoted by $\mathcal L(X_s,Y_s).$

In the case of linear functionals, i.e., when $Y=(\mathbb{R},u),$ we put $X^\flat_p=\mathcal L_{p,u}(X,\mathbb{R})$ and
$X^*=\mathcal L((X,p^s),(\mathbb{R},|\cdot|))$. The meaning of  $X^\flat_{\bar p}$ is clear.

The following proposition shows that continuous linear operators between two asymmetric normed spaces are
continuous with respect to the associated (symmetric) norm topologies.
\begin{prop}\label{p.cont-lin3}
Let $(X,p)$ and $(Y,q)$ be asymmetric seminormed spaces. Any $(p,q)$-continuous linear operator  $A:X\to Y$ is
also
$(p^s,q^s)$-continuous and the set $\mathcal{L}_{p,q}(X,Y)$ is a   cone in $\mathcal{L}(X_s,Y_s).$
Also
\begin{equation}\label{eq1.cont-lin3}
\mathcal{L}_{p,q}(X,Y)=\mathcal{L}_{\bar p,\bar q}(X,Y)\;,
\end{equation}
and
 \begin{equation}\label{eq2.cont-lin3}
 \|A|_{p^s,q^s}\le\|A|_{p,q}=\|A|_{\bar p,\bar q},\end{equation}
 for every $A\in \mathcal{L}_{p,q}(X,Y).$

In particular, every $(p,u)$-continuous linear functional is $(p^s,|\cdot|)$-continuous and $X^\flat_p$ is a
cone in the space  $X^*  $ of all continuous linear functionals on the normed space $(X,p^s).$
\end{prop}

Let us recall, following \cite{aleg-fer-greg99a}  (see also \cite[\S 1.2.4]{Cobzas}),  some notions related to Baire category in  quasi-pseudometric spaces, regarded as  bitopological spaces.
\begin{defi}\label{def.Baire}  Let $(X,d)$ be a quasi-pseudometric space. A subset $Y$ of $X$ is called:
\begin{itemize}\item
   $(d,\bar{d})$-\emph{nowhere dense} if $\,\inter_ d(\clos_{\bar{d}}(S))=\emptyset$;
\item of $(d,\bar{d})$-\emph{first category} if it is the union of a countable family of
$(d,\bar{d})$-nowhere dense sets;
\item of $(d,\bar{d})$-\emph{second category} if it is not of $(d,\bar{d})$-first category;
\item  $(d,\bar{d})$-\emph{residual} if  $X\setminus Y$ is of  $(d,\bar{d})$-first category.
\end{itemize}

 The space $(X,d)$ is called:
\begin{itemize}
 \item   $(d,\bar{d})$-\emph{Baire} if each nonempty $d$-open subset
 of $X$ is of $(d,\bar{d})$-second category;
\item  \emph{pairwise Baire} if it is both $(d,\bar{d})$-Baire  and  $(\bar{d},d)$-Baire.
\end{itemize}\end{defi}

Some versions of the uniform boundedness principle hold in the asymmetric case too.
A family $\{T_i:i\in I\}$ of continuous linear operators between two asymmetric normed spaces $(X,p)$ and $(Y,q)$ is called \emph{pointwise bounded}  if
\begin{equation*}
\sup_{i\in I} q(T_ix) < \infty \quad\mbox{for every }\; x\in X.\end{equation*}

\begin{theo}\label{t.UBP} Let $(X,p)$ and    $(Y,q)$   be asymmetric normed spaces and let  $\{T_i:i\in I\}\subseteq \mathcal{L}_{p,q}(X,Y)$ be a pointwise bounded family of continuous linear operators.
\begin{enumerate}\item[\rm 1.]{\rm(\cite[Section 2.3]{Cobzas})}\,
If  $(X,p)$ is right $p$-$K$-complete, then
\begin{equation}\label{eq2.UBP}
  \quad \sup_{i\in I}\sup_{x\in B_{p}}\bar q(T_ix) <\infty \,.\end{equation}
\item[\rm 2.] {\rm(\cite{aleg-romag12})}\,  If $(X,p)$ is of second $(p,\bar p)$-Baire category, then
\begin{equation}\label{eq3.UBP}
\sup_{i\in I}\sup_ {x\in B_{p}} q(T_ix)<\infty \,.\end{equation}
\end{enumerate}\end{theo}

The proof  (based on a lemma of Zabrejko) of the results from  Theorem \ref{t.UBP} was given in \cite{mabula-cob15}.
\begin{remark}
  Condition \eqref{eq2.UBP}  is equivalent to
  \begin{equation*}
\sup_{i\in I}\sup_ {x\in B_{\bar{p}}} q(T_ix)<\infty \,.\end{equation*}
\end{remark}

Indeed,
\begin{align*}
\sup_{i\in I}\sup_ {\bar p(x)\le 1} q(T_ix)&=\sup_{i\in I}\sup_ {p(-x)\le 1}  q(T_ix)\\
&=\sup_{i\in I}\sup_ {p(x')\le 1}  q(T_i(-x'))\\
&=\sup_{i\in I}\sup_ {p(x')\le 1} \bar q(T_ix')\,.
\end{align*}
\subsection{The dual space and the adjoint operator}
We shall consider only the case of asymmetric normed spaces. Let  $(X,p)$ be an asymmetric normed space. The space $X^\flat_p=\mathcal L_{p,u}(X,\mathbb{R})$ is called the \emph{dual  space} of $(X,p)$.
It is a cancellative cone contained in the dual of $X^*=\mathcal L_{p^s,|\cdot|}(X,\mathbb{R})$ of the normed space $(X,p^s)$. It is normed by
\begin{equation}\label{def.norm-lin-fc}\begin{aligned}
  p^\flat(\varphi)&=\sup\{u(\varphi(x)): x\in X,\, p(x)\le 1\}\\
  &=\sup\{\varphi(x): x\in X,\, p(x)\le 1\}\,,
\end{aligned}\end{equation}
for $\varphi\in X^\flat_p$.
Also
$$
\|\varphi\|^*\le p^\flat(\varphi)\,,$$
for all $\varphi\in X^*$, where
\begin{align*}
  \|\varphi\|^*&=\sup\{|\varphi(x)|:x\in X,\, p^s(x)\le 1\}\\
  &=\sup\{\varphi(x):x\in X,\, p^s(x)\le 1\}\,.
\end{align*}

In this case one can also consider the space $X^\flat_{\bar p}$ corresponding to the conjugate asymmetric norm $\bar p(x)=p(-x)$.

 On the dual $X^\flat_p$ one can define the analog of the $w^*$-topology on the dual of a normed space.
This   topology, denoted by $w^\flat$, is defined in the following way.
For $ \varphi_0\in X^\flat_p$, $x\in X$ and  $\varepsilon>0$  put

$$
V_{x,\varepsilon}(\varphi_0) =\{\varphi\in X^\flat_p: \varphi(x)-\varphi_0(x)<\varepsilon\}\,.$$

Then $w^\flat$ is the topology for which
$$
V_{x,\varepsilon}(\varphi_0),\; x\in X,\, \varepsilon>0\,,$$
is a neighborhood subbasis at $\varphi_0$, for every $\varphi_0\in X^\flat_p$.
\begin{prop}\label{p1.net-w-flat}
Let $(X,p)$ be an asymmetric normed space and let $(\varphi_i;i\in I)$  be a net in the dual $X^\flat_p$ of $(X,p)$. Then the followings are equivalent:
\begin{itemize}
\item[\rm (i)] the net $(\varphi_i;i\in I)$ is $w^\flat$-convergent to $\varphi\in X^\flat_p$;
\item[\rm (ii)] for every $x\in X$ the net $(\varphi_i(x);i\in I)$ is $u$-convergent in $\mathbb{R}$ to $\varphi(x)$;
\item[\rm (iii)] for every $x\in X$ the net $(\varphi_i(x);i\in I)$ is $|\cdot|$-convergent in $\mathbb{R}$ to $\varphi(x)$.
\end{itemize}\end{prop}

The analog of the Alaoglu-Bourbaki compactness theorem holds in this case too.
\begin{theo}\label{t.Alaog} The unit ball,
$$
B_{p^\flat}=\{\varphi\in X^\flat _p: p^\flat(\varphi)\le 1\},$$
is compact with respect to the $w^\flat$-topology.
\end{theo}

The adjoint of a continuous linear operator can be also defined within this context. Let $(X,p),(Y,q)$ be asymmetric normed spaces and  let $A\in\mathcal L_{p,q}(X,Y)$  be  a continuous linear operator.
Define $A^\flat:Y^\flat_q\to X^\flat_p$ by
$$
(A^\flat\psi)(x)=\psi(Ax)\;\mbox{ for }\; x\in X\;\mbox{ and }\; \psi\in Y^\flat_q.$$

\begin{prop}\label{p1.conj-op} Let $(X,p), (Y,q)$ be asymmetric normed spaces,  let  $A:(X,p)\to(Y,q)$  be a continuous linear operator and  let $A^\flat:Y^\flat_q\to X^\flat_p$  be  its adjoint operator.
Then $A^\flat$ is well defined,  linear and continuous   from $(Y^\flat_q,q^\flat)$ to  $(X^\flat_p,p^\flat)$ and its norm satisfies the equality
$$
\|A^\flat|_{q^\flat,p^\flat}= \|A|_{p,q}\,.$$

It is also continuous with respect to the $w^\flat$-topologies of the spaces $Y^\flat_q$ and $X^\flat_p$.
\end{prop}

A linear operator $A$  between two normed cones $(X,p)$ and $Y,q)$ is called compact if the set $A(B_p)$ is relatively compact in $Y$.
The following is an asymmetric analog of Schauder's compactness theorem.
\begin{theo}[see \cite{cobz06} or \cite{Cobzas}, \S 2.4.4]\label{t.Schaud-asym}  Let $(X,p), (Y,q)$ be  asymmetric normed spaces, let $A:(X,p)\to(Y,q)$ be a continuous linear operator and let $A^\flat:Y^\flat\to X^\flat$ be its adjoint. If $A$ is compact, then $A^\flat$ is also compact.
\end{theo}

\section{Bilinear operators on asymmetric normed spaces}
We introduce in this section continuous bilinear operators on asymmetric normed spaces and the corresponding spaces of continuous bilinear operators.
\subsection{Bilinear operators}
Let $X,Y,Z$ be real vector spaces. A mapping $T:X\times Y\to Z$ is called a \emph{bilinear operator} if for every $y\in Y,\; T(\cdot,y):X\to Z$ is a linear
operator and for every $x\in X,\; T(x,\cdot):Y\to Z$ is also a linear operator. If  $(X,p_1),(Y,p_2)$ and $(Z,q)$ are asymmetric normed spaces, then we are interested in the continuity
properties of a bilinear operator $T:X\times Y\to Z$. We start with the following simple result.
\begin{lemma}\label{le1.asyn-bilin} Let $\beta\ge 0$ and $r>0$. The  inequalities
\begin{equation}\label{eq1.asyn-bilin}\begin{aligned}
  {\rm (i)}\quad &q(T(x,y))\le\frac{\beta}{r^2}\,p_1(x)\,p_2(y)\;\mbox{for all }\; (x,y)\in X\times Y, \;\mbox{and}\\
  {\rm (ii)}\quad &q(T(x,y))\le\beta\;\mbox{ for all }\; (x,y)\in X\times Y \;\mbox{with }\; p_1(x)\le r\;\mbox{ and }\; p_2(y)\le r,
\end{aligned}\end{equation}
are equivalent.
\end{lemma}\begin{pf}
  (i)\;$\Rightarrow$\; (ii).

  If $p_1(x)\le r$ and $p_2(y)\le r$, then, by (i),
  $$
  q(T(x,y))\le\frac{\beta}{r^2}\,p_1(x)p_2(y)\le\frac{\beta}{r^2}\, r^2=\beta\,.$$

  (ii)\;$\Rightarrow$\; (i).

 Let $(x,y)\in X\times Y$. If $p_1(x)p_2(y)>0$, then
  $$
  x'=\frac{r}{p_1(x)}  x\quad\mbox{and}\quad  y'=\frac{r}{p_2(y)}  y\,,$$
  satisfy $p_1(x')=r=p_2(y')$, so that
  $$
  \frac{r}{p_1(x)} \, \frac{r}{p_2(y)}\, q(T(x,y))=q(T(x',y'))\le \beta,\,,$$
   which implies
  $$q(T(x,y))\le\frac{\beta}{r^2}\,p_1(x)p_2(y)\,.$$

  If, for instance, $p_1(x)=0\le r $ and $p_2(y)>0$, then $p_1(tx)=0$ for all $t>0$, so that, taking $y'$ as above,
  $$
  \frac{r}{p_2(y)}\,t\, q(T(x,y))=q(T(tx,y'))\le \beta\,.$$

  Hence
  $$
  0\le q(T(x,y))\le\frac{\beta p_2(y)}{tr}\,\to 0\;\mbox{ as }\; t\to\infty\,.$$

  This shows that $$q(T(x,y))=0=\frac{\beta}{r^2}\, p_1(x) p_2(y)\,.$$

  If $p_1(x)=0=p_2(y)$, then one can work with $x'=tx,\,  y'=ty$, for $t>0$, to deduce that $q(T(x,y))=0.$
\qed\end{pf}

\begin{prop}\label{p2.asyn-bilin} Let $(X,p_1),(Y,p_2)$, $(Z,q)$ be asymmetric normed spaces and  let $T:X\times Y\to Z$  be a bilinear operator.
  The followings are equivalent:
\begin{itemize}
\item[\rm (i)] the mapping $T$ is continuous;
\item[\rm (ii)] the mapping $T$ is continuous at $(0,0)$;
\item[\rm (iii)]  there exists $\beta\ge 0$ such that
\begin{equation}\label{eq1.bilin-semi-Lip}
q(T(x,y))\le \beta p_1(x)p_2(y)\;\mbox{ for all }\; (x,y)\in X\times Y.\end{equation}
\end{itemize}\end{prop}\begin{pf}
  (i)\;$\Rightarrow$\;(ii). \, This is obvious.

(ii)\;$\Rightarrow$\;(iii).
  Using the continuity of $T$ at $(0,0)$ it follows that there exists $r>0$ such that
  $$
  q(T(x,y))\le 1\;\mbox{ for all }\; (x,y)\in rB_{p_1}\times rB_{p_2}\,.$$

  But then, by Lemma \ref{le1.asyn-bilin},
  $$
  q(T(x,y))\le \frac1{r^2}p_1(x)p_2(y)\,,$$
for all $ (x,y)\in X\times Y$.

  (iii)\;$\Rightarrow$\;(i).

  Let $(x,y)\in X\times Y$. We show that $T(x_n,y_n)\xrightarrow{q}T(x,y)$ for every  sequence $(x_n,y_n)\in X\times Y,\, n\in\mathbb{N},$  which is $p_1\times p_2$-convergent to $(x,y)$.
  This is equivalent to $$p_1(x_n-x)\to 0\quad\mbox{and}\quad p_2(y_n-y)\to 0. $$ The inequality
  $$p_2(y_n)\le p_2(y_n-y)+p_2(y),\; n\in\mathbb{N}\,,$$
  implies the boundedness of the sequence $(p_2(y_n))_{n\in\mathbb{N}}\,.$

  Writing
  $$
  T(x_n,y_n)-T(x,y)=T(x_n-x,y_n)+T(x,y_n-y)\,,$$
  it follows  that
\begin{align*}
q(T(x_n,y_n)-T(x,y))&\le q(T(x_n-x,y_n))+q(T(x,y_n-y))\\
&\le\beta (p_1(x_n-x)p_2(y_n)+p_1(x)p_2(y_n-y))\longrightarrow 0 \,,
\end{align*}
  as $\, n\to\infty$. \qed\end{pf}

\begin{remark}
  As in the case of linear operators (see Section \ref{SS.lin-ops-asyn}), a number $\beta\ge 0$ for which \eqref{eq1.bilin-semi-Lip} holds is called a \emph{semi-Lipschitz constant} for $T$.
\end{remark}

For a bilinear operator $T:X\times Y\to Z$ put
\begin{equation}\label{def.norm-bilin}
\|T|=\|T|_{p_1,p_2;q}=\sup\{q(T(x,y)): (x,y)\in B_{p_1}\times B_{p_2}\}\,.\end{equation}
Denote by  $\L^2_{p_1,p_2;q}(X,Y;Z)$ the set of all continuous bilinear operators from $(X\times Y,p_1\times p_2)$ to $(Z,q)$.

Denote by   $\mathcal L^2(X,Y;Z)$   the space $\L^2_{p^s_1,p^s_2;q^s}(X,Y;Z)$ and   let $\|T\|$  be the norm of a bilinear operator in  $\mathcal L^2(X,Y;Z)$  given by
\begin{equation}
 \label{def.norm-bilin-sym}\|T\|=\sup\{q^s(T(x,y)): (x,y)\in X\times Y,\, p_1^s(x)\le 1,\, p_2^s(y)\le 1\}\,.
\end{equation}

Then $\mathcal L^2(X,Y;Z)$ is a linear space  normed by \eqref{def.norm-bilin-sym}, which is Banach if  $(Z,q)$ is a bi-Banach space (i.e.  if $(Z,q^s)$ is complete).

\begin{prop}\label{p3.asyn-bilin}  Let $(X,p_1),(Y,p_2)$, $(Z,q)$ be asymmetric normed spaces and  let $T:X\times Y\to Z$  be a bilinear operator.
\begin{enumerate}\item[\rm 1.]  The operator $T$ is continuous if and only if $\,\|T|_{p_1,p_2;q}<\infty$. In this case  $\|T|_{p_1,p_2;q}$ is the smallest semi-Lipschitz constant for $T$.
\item[\rm 2.]
The space $\L^2_{p_1,p_2;q}(X,Y;Z)$ is a   cone   in the space  $\mathcal L^2(X,Y;Z)$. The functional   $\|\cdot|_{p_1,p_2;q}$ is a norm on it and
$$ \|T\|\le\|T|_{p_1,p_2;q}\,,$$
for all  $T\in\L^2_{p_1,p_2;q}(X,Y;Z)$.
\end{enumerate}\end{prop}\begin{pf} 1. For convenience, we denote the norm of $T$ by $\|T|$. If $\|T|<\infty$, then by Lemma \ref{le1.asyn-bilin},
$$q(T(x,y))\le\|T|\,p_1(x)p_2(y)\;\mbox{ for all }\;(x,y)\in X\times Y.$$

By Proposition \ref{p2.asyn-bilin}, this  shows that $T$ is continuous and that $\|T|$ is a semi-Lipschitz constant for $T$.
If $\beta$ is a semi-Lipschitz constant for $T$, then, appealing again to Lemma \ref{le1.asyn-bilin} (for $r=1$), it follows $\|T|\le\beta$.\smallskip

2.\;
  Let $(x,y)\in X\times Y$. Then
  $$
  q(T(x,y))\le \|T|\,p_1(x)p_2(y)\le  \|T|\,p^s_1(x)p^s_2(y)\,,$$
  and
  $$
  \bar q(T(x,y))=q(T(-x,y))\le \|T|\,p_1(-x)p_2(y)\le  \|T|\,p^s_1(x)p^s_2(y)\,,$$
  implying
  $$
  q^s(T(x,y)) \le  \|T|\,p^s_1(x)p^s_2(y)\,.$$

  Consequently, $T\in\mathcal L^2(X,Y;Z)$ and $\|T\|\le\|T|_{p_1,p_2;q}$.
  The fact that $\L^2_{p_1,p_2;q}(X,Y;Z)$ is a cancellative cone in $\mathcal L^2(X,Y;Z)$ is easily seen.
\qed\end{pf}

It is known that in the case of Banach spaces a separately continuous multilinear   operator is globally continuous. This result was extended to the asymmetric case by Latreche and Dahia \cite{latreche20}. For  the sake of completeness, we present here the simple proof in the case of bilinear operators.

\begin{prop}\label{p1.bilin-sep-cont} Let  $(X,p_1)$ be  an asymmetric normed space of second $(p_1,\bar p_1)$-Baire category and  let $(Y,p_2)$, $(Z,q)$  be  asymmetric normed spaces. If $T:X\times Y\to Z$ is a separately continuous bilinear operator, then
$T$ is  $(p_1\times p_2,q)$-continuous.
\end{prop}\begin{pf}
  Suppose that $T$ is not continuous at $(0,0)\in X\times Y$. Then there exist $r>0$ and a sequence $\big((x_n,y_n)\big)_{n\in\mathbb{N}}$\,that is $p_1\times p_2$-convergent to (0,0),
such that\begin{equation}\label{eq1.bilin-sep-cont}
q(T(x_n,y_n))\ge r\quad\mbox{for all }\; n\in\mathbb{N}.\end{equation}

Let $A_n:X\to Z$ be defined by $A_n(\cdot)=T(\cdot,y_n),\, n\in\mathbb{N}$. Then $A_n\in\mathcal L_{p_1,q}(X,Z)$ for every $n\in\mathbb{N}$, so there exists $\lambda_n\ge 0$ such that
\begin{equation}\label{eq1a.bilin-sep-cont}
q(A_nx)\le \lambda_np_1(x) \quad\mbox{for all }\; x\in X\,.\end{equation}

For every $x\in X$,
$$A_nx=T(x,y_n)\xrightarrow{q} 0\;\mbox{ as }\; n\to\infty\,,$$
so that
$$\sup_nq(A_nx)<\infty  \quad\mbox{for all }\; x\in X\,,$$
that is, the family $\{A_n:n\in\mathbb{N}\}$ is pointwise bounded. By Theorem \ref{t.UBP}.2
$$
\beta:=\sup_n\sup_{x\in B_{p_1}}q(T(x,y_n))<\infty\,.$$

By \eqref{eq1a.bilin-sep-cont},
$$
0<r\le q(T(x_n,y_n))\le \lambda_np_1(x_n)\,,$$
implying $p_1(x_n)>0$ for all $n$. But  then,
$$
q\left(T\left(x_n/p_1(x_n),y_n\right)\right)\le \beta\,,$$
  for all $n$.

 Hence, by \eqref{eq1.bilin-sep-cont},
 $$
 0<r\le  q(T(x_n,y_n))\le\beta p_1(x_n)\to 0\;\mbox{ as }\; n\to\infty\,,$$
 a contradiction which shows that $T$ must be continuous  at $(0,0)$ and so, continuous on $X\times Y$.
\qed\end{pf}
\begin{remark} Proposition \ref{p1.bilin-sep-cont} is valid when $(X,p_1)$ is a normed space of second Baire category, in particular when it  is a Banach space.
\end{remark}
\subsection{Bilinear forms and the $w^2$-topology}

Let $u$ be  the asymmetric norm on $\mathbb{R}$  introduced in  Example \ref{ex1.Ru}.
Consider now the case of bilinear forms, i.e., the space $\L^2_{p_1,p_2;u}(X,Y)$ of  bilinear forms from $X\times Y$ to $\mathbb{R}$
which are $p_1\times p_2$ to $u$ continuous.  For the sake of simplicity we shall denote it by  $\L^2(X,Y).$ It is a cancellative cone, normed by
\begin{equation}\label{def.norm-bilin-fcs}\begin{aligned}
  \|b|_{p_1,p_2}&=\sup\{u(b(x,y)) : (x,y)\in B_{p_1}\times B_{p_2}\}\\
  &=\sup\{b(x,y) : (x,y)\in B_{p_1}\times B_{p_2}\}\,,
\end{aligned}\end{equation}
for $b\in\L^2(X,Y).$
\begin{remark}\label{re.as-bilin-fcs}
  By Proposition \ref{p3.asyn-bilin},  $\,\L^2(X,Y)=\L^2_{p_1,p_2;u}(X,Y)$ is contained in the space $\mathcal L^2(X,Y)=\L^2_{p_1^s,p_2^s;|\cdot|}(X,Y;\mathbb{R})$ and
  \begin{equation}\label{eq.as-bilin-fcs}
  \|b\|\le \|b|_{p_1,p_2}\,,\end{equation}
  where $\|b\|=\|b|_{p^s_1,p^s_2}$.
\end{remark}

The functional $\delta:\L^2(X,Y)\times\L^2(X,Y)\to[0,\infty]$ defined for $b_1,b_2\in \L^2(X,Y)$ by
$$
\delta(b_1,b_2)=\begin{cases} \|b_2-b_1|_{p_1,p_2}\quad&\mbox{if }\; b_2\in b_1+ \L^2(X,Y)\\
\infty &\mbox{otherwise}
\end{cases}$$
is an extended quasi-pseudometric on the cone $\L^2(X,Y).$

Now we shall introduce on  $\L^2(X,Y)$ the analog of the $w^\flat$ topology on the dual of an asymmetric normed space.
Denote by $\theta$ the null bilinear form. For $(x,y)\in X\times Y$, $\varepsilon>0$ and $b_0\in \L^2(X,Y)$ put
$$  V_{x,y;\varepsilon}(b_0) =\{b\in \L^2(X,Y): b(x,y)-b_0(x,y)<\varepsilon\}=b_0+V_{x,y;\varepsilon}(\theta) \,.
$$

Denote by $w^2$ the topology on $  \L^2(X,Y)$ for which the family of sets
$$
V_{x,y;\varepsilon}(b_0),\;\; (x,y)\in X\times Y,\, \varepsilon>0\,,$$
is a neighborhood subbasis at each $b_0\in  \L^2(X,Y)$.

\begin{remark}\label{re.w2} Since
$$
W_{x,y;\varepsilon}(b_0):=V_{x,y;\varepsilon}(b_0)\cap V_{-x,y;\varepsilon}(b_0)=\{b\in \L^2(X,Y): |b(x,y)-b_0(x,y)|<\varepsilon\}\,,$$
it follows that the family of sets$$W_{x,y;\varepsilon}(b_0),\; \varepsilon>0,\; (x,y)\in X\times Y\,,$$
is also a subbasis of $w^2$-neighborhoods of $b_0$.
\end{remark}

The $w^2$-convergence of a net in $\L^2(X,Y)$     can be characterized in the following way.

\begin{prop}\label{p1.w2-conv}  Let $(b_i; i\in I)$ be a net   in $\L^2(X,Y)$ and let $b\in\L^2(X,Y)$. Then the followings are equivalent:
\begin{itemize}\item[\rm (i)]  the net $(b_i; i\in I)$ is $w^2$-convergent to $b$;
\item[\rm (ii)] for every $(x,y)\in X\times Y,\; b_i(x,y)\overset{u}{\longrightarrow} b(x,y)$\;in\ $\mathbb{R}$;
\item[\rm (iii)]for every $(x,y)\in X\times Y,\;  b_i(x,y) \overset{|\cdot|}{\longrightarrow} b(x,y)$\;in\ $\mathbb{R}$.
\end{itemize}\end{prop}
\begin{pf}  Notice that, for every $\varepsilon>0$,
$$
\alpha<\varepsilon\iff u(\alpha)<\varepsilon\,$$
for all $\alpha\in\mathbb{R}.$

 (i)$\iff$(ii).

This follows from the equivalences
  \begin{align*}
    b_i\overset{w^2}{\longrightarrow}b&\iff \forall (x,y)\in X\times Y,\; \forall \varepsilon>0,\;\exists i_0\in I,\;\forall i\ge i_0,\;\; b_i\in V_{x,y:\varepsilon}(b)\\
    &\iff \forall (x,y)\in X\times Y,\; \forall \varepsilon>0,\;\exists i_0\in I,\;\forall i\ge i_0,\;\; b_i(x,y)-b(x,y)< \varepsilon\\
    &\iff \forall (x,y)\in X\times Y,\;  b_i(x,y)\overset{u}{\longrightarrow}b(x,y)\,.
  \end{align*}

(ii)\;$\iff$\; (iii).

The implication (iii)\;$\Rightarrow$\;(ii) follows from the inequality
$$
u(b_i(x,y)-b(x,y))\le |b_i(x,y)-b(x,y)|\,.$$

To prove the reverse implication,
suppose now that $b_i(x,y)\overset{u}{\longrightarrow}b(x,y)$ for all $(x,y)\in X\times Y$.

Let $(x,y)\in X\times Y$ and $\varepsilon>0$.  Then there exists $i_0\in I$ such that
\begin{align*}
 &b_i(x,y)-b(x,y)<\varepsilon\;\mbox{ and }\\  -&b_i(x,y)+b(x,y)=b_i(-x,y)-b(-x,y)<\varepsilon\,,
\end{align*}
for all $i\ge i_0$, which implies that
$$
|b_i(x,y)-b(x,y)|<\varepsilon\,,$$
for all $i\ge i_0$, that is, $b_i(x,y) \overset{|\cdot|}{\longrightarrow}b(x,y)$
\qed\end{pf}

The analog of the Alaoglu-Bourbaki theorem holds in this case too.
\begin{theo}\label{t.Alaog-bilin} The closed unit ball
$$
B=\{b\in \L^2_{p_1,p_2;u}(X,Y): \|b|_{p_1,p_2;u}\le 1\}$$
is $w^2$-compact.
\end{theo}\begin{pf} For convenience, denote the space $\L^2_{p_1,p_2;u}(X,Y)$ by $\L^2$ and the corresponding norm by $\|\cdot|$, i.e.
$\|b|=\|b|_{p_1,p_2;u}$ for  $b\in\L^2$.

Consider the space $\Lambda=\mathbb{R}^{X\times Y}$ endowed with the product topology, denoted by $\tau$.
The convergence of a net $(\lambda_i:i\in I)$ in $\Lambda$ to $\lambda\in\Lambda$ can be characterized by:
$$
\lambda_i\overset{\tau}{\longrightarrow}\lambda\iff \forall (x,y)\in X\times Y,\; \lambda_i(x,y)\overset{|\cdot|}{\longrightarrow} \lambda(x,y)\;\mbox{ in }\ \mathbb{R}\,.$$

But then, the equivalence (i)$\iff$(iii) from Proposition \ref{p1.w2-conv} shows that
\begin{equation}\label{eq.tau-w2}
  b_i\overset{w^2}{\longrightarrow}b\iff b_i\overset{\tau}{\longrightarrow}b\,,
\end{equation}
for every net $(b_i)$ in $\L^2_{p_1,p_2;u}(X,Y)$ and every $b\in \L^2_{p_1,p_2;|\cdot|}(X,Y).$

Let $\Delta_{x,y}=[-p^s_1(x)p^s_2(y),p^s_1(x)p^s_2(y)]$, $ (x,y)\in X\times Y,$ and  let $\Gamma=\prod_{(x,y)}\Delta_{x,y}$.  Since each $\Delta_{x,y}$ is compact in $(\mathbb{R},|\cdot|)$ it follows that $\Gamma$  is compact with respect to the product topology, that is, $\tau$-compact.

The inequalities \begin{align*}
  &b(x,y)\le p_1^s(x) p_2^s(y)\;\mbox{ and }\\-&b(x,y)=b(-x,y)\le p_2^s(y) p_1^s(x)
\end{align*}
show that $|b(x,y)|\le p_1^s(x) p_1^s(x)$, that is,
$b(x,y)\in\Delta_{x,y}$ for all $(x,y)\in X\times Y$, or, equivalently
$$B\subseteq\Gamma\,.$$

\emph{Fact} I. $B$ \emph{is a} $\tau$-\emph{closed subset of} $\Gamma$. \smallskip

If $(b_i:i\in I)$ is a net   in $B$ that is $\tau$-convergent to $\gamma\in\Gamma$, then
 $$
b_i\overset{\tau}{\longrightarrow}\gamma\iff \forall (x,y)\in X\times Y,\; b_i(x,y)\overset{|\cdot|}{\longrightarrow}\gamma(x,y)\,.$$

For every  $ x,x'\in X,\, y\in Y$ and $s,t\in\mathbb{R}$ the equality
$$
b_i(sx+tx',y)=sb_i(x,y)+tb_i(x',y)\,,$$
holds for every $i\in I$. Passing to limit for  $i\in I$, one obtains
$$
\gamma(sx+tx',y)=s\gamma(x,y)+t\gamma(x',y)\,,$$
  for all $x,x'\in X,\, y\in Y$.
  The linearity  with respect to $y$ can be proved analogously, so that    $\gamma $ is a bilinear form on $X\times Y$.

  Similarly, the inequality
  $$
  b_i(x,y)\le p_1(x)p_2(y)\ $$
  yields
  $$
  \gamma(x,y)\le p_1(x)p_2(y)\,, $$
  for all $(x,y)\in X\times Y$, showing that $\gamma \in B$.

 As a $\tau$-closed subset of the $\tau$-compact set $\Gamma$,   the ball $B$ is  $\tau$-compact  too.
 By \eqref{eq.tau-w2},  the topologies $\tau$ and $w^2$ agree on $B$, so that $B$ is $w^2$-compact.
\qed\end{pf}

\subsection{The adjoint of a bilinear operator}
Let $(X,p_1), (Y,p_2)$  and $(Z,q)$ be asymmetric normed spaces and  let $A\in\L_{p_1,p_2;q}^2(X,Y;Z)$.  The adjoint of the operator $A$ is the operator $A^\flat:Z^\flat_q\to \L^2(X,Y)$
defined for $\psi \in Z^\flat_q$ by
\begin{equation}\label{def.conj-bilin}
(A^\flat\psi)(x,y)=\psi(A(x,y)),\; (x,y)\in X\times Y\,.\end{equation}
\begin{prop}\label{p1.conj-bilin} Let $(X,p_1), (Y,p_2)$, $(Z,q)$ and   let  $ A$ be as above.
\begin{enumerate}\item[\rm 1.] The operator $A^\flat$ is well defined,  linear (meaning additive and positively homogeneous) and continuous from $(Z^\flat_q,q^\flat)$ to $(\L^2(X,Y),\|\cdot|)$ and
$$
\|A^\flat|_{q^\flat,\|\cdot|}=\|A|_{p_1,p_2;q}\,.$$
  \item[\rm 2.] The operator $A^\flat$ is also continuous with respect to the topologies $w^\flat$ on $Z^\flat_q$ and $w^2$ on $\L^2(X,Y)$.
  \end{enumerate}\end{prop}\begin{pf}
  1.\; For simplicity we shall denote the norms by $\|\cdot|$, without subscripts.

The linearity of $A^\flat$ is obvious.

  For $\psi\in Z^\flat_q$ and $(x,y)\in X\times Y$ we have
 \begin{align*}
     (A^\flat\psi)(x,y)&=\psi(A(x,y))\le q^\flat(\psi)q(A(x,y))\\
     &\le  q^\flat(\psi)\|A|p_1(x)p_2(y)\,.
 \end{align*}

By Proposition \ref{p2.lin-op-cone} this implies
 $$
 \|A^\flat\psi|\le \|A| q^\flat(\psi)$$
 so that $A^\flat\in \mathcal L_{q^\flat,\|\cdot|}(Z^\flat,\L^2(X,Y))$ and
 $$
 \|A^\flat|\le \|A| \,.$$

 By Theorem 2.2.2 in \cite{Cobzas}, for any $(x,y)\in B_{p_1}\times B_{p_2}$ with $A(x,y)\ne 0$, there exists $\psi\in Z^\flat$ such that
 $$
 q^\flat(\psi)=1\quad\mbox{and}\quad \psi(A(x,y))=q(A(x,y))\,.$$

 But then
 $$
 q(A(x,y))=\psi(A(x,y))=(A^\flat\psi)(x,y)\le \|A^\flat|p_1(x)p_2(y)\le\|A^\flat|\,.$$

 Passing to supremum with respect to $(x,y)\in B_{p_1}\times B_{p_2}$ one obtains
 $$ \|A|\le\|A^\flat|\,.$$

  2.\; We have to show that
  $$
  \psi_i\overset{w^\flat}{\longrightarrow} \psi \;\Rightarrow\; A^\flat\psi_i\overset{w^2}{\longrightarrow} A^\flat\psi\,,$$
  for every net $(\psi_i:i\in I)$ in $Z^\flat_q$ which is $w^\flat$-convergent to $\psi \in Z^\flat_q$.

  Taking into account Proposition \ref{p1.w2-conv}, this follows from the following relations
  \begin{align*}
     \psi_i\overset{w^\flat}{\longrightarrow} \psi &\iff \forall z\in Z,\; \psi_i(z)\overset{u}{\longrightarrow} \psi(z)\\
     &\quad \Rightarrow\; \forall (x,y)\in X\times Y,\; \psi_i(A(x,y))\overset{u}{\longrightarrow} \psi(A(x,y))\\
     &\iff \forall (x,y)\in X\times Y,\; (A^\flat\psi_i)(x,y)\overset{u}{\longrightarrow} (A^\flat\psi)(x,y)\\
     &\iff  A^\flat\psi_i\overset{w^2}{\longrightarrow} A^\flat\psi\,.
  \end{align*}\qed\end{pf}
\section{Precompact bilinear operators}
In this section we introduce the precompactness notions for bilinear operators on asymmetric normed spaces and prove some of their fundamental properties, including a Schauder-type theorem.
\subsection{Definition and some properties}

According   to the precompactness notions from Definition \ref{def.prec} one can consider    the corresponding precompactness notions for  bilinear operators.
\begin{defi}\label{def.prec-op}
 Let $(X,p_1),(Y,p_2)$ and   $(Z,q)$ be  quasi-pseudometric spaces. A bilinear operator $T:X\times Y\to Z$ is called:
  \begin{itemize}
  \item $q$-\emph{precompact} if $T(B_{p_1}\times B_{p_2})$ is a $q$-precompact subset of $Z$;
  \item $q^s$-\emph{precompact} if $T(B_{p_1}\times B_{p_2})$ is a $q^s$-precompact subset of $Z$;
  \item \emph{h-precompact} if $T(B_{p_1}\times B_{p_2})$ is a hereditarily $q$-precompact subset of $Z$.
  \end{itemize}\end{defi}

  Here $q$-precompactness is understood as precompactness with respect to the quasi-metric $d_q(x,y)=q(y-x)$ induced by the asymmetric norm $q$, and similarly for $q^s$.
Recall that we denote by $B_p$ the closed unit ball of an asymmetric normed space $(W,p)$.

  If $q$ is an asymmetric norm on a cone $Z$, then by $q$-precompactness one understands the precompactness with respect to the extended quasi-metric    \eqref{def.ext-dist}
  corresponding to $q$.

We denote by  $ \mathcal K^2_{p_1,p_2;q}(X,Y;Z)\, (\mathcal K^2_{p_1,p_2;q^s}(X,Y;Z))$  the set of all  $q$-precompact  ($q^s$-precompact) bilinear operators and by $ h\mathcal K^2_{p_1,p_2;q}(X,Y;Z)$  the set of all  $h$-precompact bilinear operators.

The following remark shows that precompact bilinear operators are actually continuous.
\begin{remark}
  The followings hold:\begin{itemize}
  \item[\rm (i)]
  a $q^s$-precompact bilinear operator is $q$-precompact;
  \item[\rm (ii)] a $q$-precompact bilinear operator is $(p_1\times p_2,q)$-continuous;
\item[\rm (iii)] a $q^s$-precompact bilinear operator is $(p_1\times p_2,q^s)$-continuous.
   \end{itemize}
\end{remark}
The assertion (i) follows from Remark \ref{re.prec}.

The assertion (ii) and (iii) follow from Remark \ref{re.prec-bd}. Indeed,

$$
\|T|_{p_1,p_2;q}=\sup q\left(T(B_{p_1}\times B_{p_2})\right)<\infty\,,$$
and similarly for $q^s$.

\subsection{A Schauder-type theorem}
The following result is the asymmetric bilinear  analog of the famous Schauder's  theorem on the compactness of the adjoint of a compact linear operator.
Recall that it also holds for compact linear operators on  asymmetric normed spaces (see Theorem \ref{t.Schaud-asym}).

Denote by $\|\cdot|$ the asymmetric norm  $\|\cdot|_{p_1,p_2;u}$ on $\L^2_{p_1,p_2}(X,Y)$.

\begin{theo}\label{t.Schaud-comp-bilin}
If the operator $T\in\L^2_{p_1,p_2;q}(X,Y;Z)$ is $q^s$-precompact, then the adjoint operator $T^\flat:Z^\flat_q\to\L^2_{p_1,p_2}(X,Y)$
  is   $\|\cdot|$-precompact.
\end{theo}\begin{pf}
  Let $\varepsilon>0$. By the $q^s$-precompactness of $T$ there exist $n\in\mathbb{N}$  and $$(x_1,y_1),\dots,(x_n,y_n) \in B_{p_1}\times B_{p_2}$$ such that
\begin{equation}\label{eq1.Schaud-asyn}
\forall (x,y) \in B_{p_1}\times B_{p_2},\; \exists k_0\in\{1,\dots,n\}\;\mbox{ s.t. }\; q^s(T(x,y)-T(x_{k_0},y_{k_0}))\le \varepsilon\,.\end{equation}

For  $\psi\in B_{q^\flat}$
$$
W(T^\flat\psi)=\{b\in \L^2_{p_1,p_2}(X,Y):b(T(x_k,y_k))-(T^\flat\psi)(x_k,y_k)<\varepsilon,\, k=1,\dots,n\}$$
is a $w^2$-open neighborhood of $T^\flat\psi$ and
$$
W(T^\flat\psi),\;\; \psi\in B_{q^\flat},$$
is a $w^2$-open cover of the set  $T^\flat(B_{q^\flat})$. By Alaoglu-Bourbaki theorem (Theorem \ref{t.Alaog}) the ball $B_{q^\flat}$ is $w^\flat$-compact.
By Proposition \ref{p1.conj-bilin} the operator $T^\flat$ is also $w^\flat$ to $w^2$-continuous, so that the set  $T^\flat(B_{q^\flat})$ is $w^2$-compact.
It follows that there exists $m\in\mathbb{N}$ and $\psi_1,\dots,\psi_m\in B_{q^\flat}$ so that for every $ \psi\in B^\flat_q$ there exists $i\in\{1,\dots,m\}$ such that
\begin{equation}\label{eq2.Schaud-bilin}
\;(T^\flat\psi)(x_k,y_k)-(T^\flat\psi_i)(x_k,y_k)<\varepsilon,\; k=1,\dots,n\,.
\end{equation}

Let  $\psi\in B^\flat$.    Choose first $\psi_i\in B^\flat$, for some $i \in\{1,\dots,m\},$ according to \eqref{eq2.Schaud-bilin}.

Now, for $(x,y)\in B_{p_1}\times B_{p_2}$, pick $k_0\in\{1,\dots,n\}$, according to \eqref{eq1.Schaud-asyn}.  Then
\begin{align*}
  (T^\flat\psi-T^\flat\psi_i)(x,y)&= \psi(T(x,y)-T(x_{k_0},y_{k_0}))\\
  &+ (T^\flat\psi)(x_{k_0},y_{k_0})-(T^\flat\psi_i)(x_{k_0},y_{k_0})\\
   &+ \psi_i(T(x_{k_0},y_{k_0})-T(x,y))\\
   &<\varepsilon+2q^s(T(x,y)-T(x_{k_0},y_{k_0}))\\
     &\le 3\varepsilon\,.
\end{align*}

It follows
$$
(T^\flat\psi-T^\flat\psi_i)(x,y)\le 3\varepsilon\,,$$
for all $(x,y)\in B_{p_1}\times B_{p_2}$. By Proposition \ref{p3.asyn-bilin} this implies
$$
\|T^\flat\psi-T^\flat\psi_i|_{p_1,p_2}\le 3\varepsilon\,,$$
showing that $T^\flat\psi_i,\, i=1,\dots,m$ is a $3\varepsilon$-net for $T^\flat(B_{q^\flat})$ with respect to the asymmetric norm \eqref{def.norm-bilin-fcs}  on $\L^2(X,Y).$
\qed\end{pf}

\subsection{Ideal properties} It is known that the space of compact linear operators on a Banach space $X$ is a closed ideal in the space of all continuous linear
operators on $X$.  Pietsch \cite{pietsch83,pietsch84}  defined the ideals of multilinear operators and sketched a program of work. His ideas were developed by several authors. We restrict the presentation to the case of bilinear operators.

 For each triple  $(X,Y;Z)$ of normed spaces   let $\mathcal L^2(X,Y;Z)$ and $\mathcal K^2(X,Y;Z)$ be the spaces of all continuous bilinear operators  and  of all precompact bilinear operators from $X\times Y$ to $Z$, respectively.
Denote by $\mathcal B^3$ the class of all   triples  $(X,Y;Z)$ of Banach spaces.  Consider the mapping $\mathcal L^2$ which attributes to each triple  $(X,Y;Z)\in \mathcal B^3$
the space $\mathcal L^2(X,Y;Z)$.

\begin{defi}\label{def.bid} A \emph{bideal} is a mapping $\mathcal I^2$ defined on $\mathcal B^3$ such that  the following conditions are fulfilled:
\begin{enumerate}\item[\rm (a)] $\mathcal I^2(X,Y;Z)$ is a linear subspace of $\mathcal L^2(X,Y;Z)$ for every $(X,Y;Z)\in \mathcal B^3$;
\item[\rm (b)] for every  $(X,Y;Z)\in \mathcal B^3$ and every Banach space $Z_1$,  $RT\in \mathcal I^2(X,Y;Z_1)$ for all $T\in  \mathcal I^2(X,Y;Z)$ and $R\in\mathcal L(Z,Z_1)$;
\item[\rm (c)] for all  $(X,Y;Z)\in \mathcal B^3$ and  all Banach spaces  $X_1,Y_1$,  \, $T\circ(S_1,S_2)\in \mathcal I^2(X,Y;Z)$ for all $T\in \mathcal I^2(X,Y;Z),\, S_1\in\mathcal L(X_1,Y)$ and $S_2\in \mathcal L(Y_1,Y)$, where $(S_1,S_2):X_1\times X_2\to X\times Y$ is given by
$$
(S_1,S_2)(x_1,y_1)=(S_1x_1,S_2y_1), \quad (x_1,y_1)\in X_1\times Y_1\,.$$
\end{enumerate}

A bideal $\mathcal I^2$ is called \emph{closed} if $\mathcal I^2(X,Y;Z)$ is  a closed subspace of $\mathcal L^2(X,Y;Z)$ for every $(X,Y;Z)\in\mathcal B^3$. \end{defi}

 Apparently unaware of Pietsch's work, Ramanujan and Schock \cite{schock85} also introduced a notion of bideal, where, besides  the conditions (a)--(c) from Definition \ref{def.bid}, the following one is added:
 \begin{enumerate}\item[\rm (d)]  for every $(X,Y;Z)\in \mathcal B^3$,\, $\varphi\otimes z$ belongs to $\mathcal I^2(X,Y;Z)$ for all continuous bilinear forms $\varphi\in \mathcal L^2(X,Y;\mathbb{R})$ and all $z\in Z$, where
$$
(\varphi\otimes z)(x,y)=\varphi(x,y)z,\quad \mbox{for }\;(x,y)\in X\times Y. $$
\end{enumerate}

Ramanujan and Schock \cite{schock85} proved the following result.
\begin{theo} The compact bilinear operators form a closed bideal in $\mathcal L^2$.\end{theo}

The paper mentioned above  also contains  other examples of bideals: bilinear operators with finite range and the bilinear analogs of nuclear and absolutely summing linear operators. A factorization result is also proved: a bilinear operator is precompact if and only if it admits a factorization through some subspaces of $c_0$.
The proof of this last result contained some  flaws which were fixed  by Ruch \cite{ruch89} (see also \cite{ruch96}).  Recall that a   bilinear operator $T:X\times Y\to Z$ is called compact if $T(B_X\times B_Y)$ is relatively compact in $Z$. The idea to use precompact bilinear operators, i.e. bilinear operators $T$ for which $T(B_X\times B_Y)$ is precompact in the normed space $Z$, belongs to Ruch \cite{ruch89}, who also gave an example of a precompact bilinear operator which is not compact. Obviously, if the address space $Z$ is Banach then compact and precompact bilinear operators agree.

In what follows we shall examine the closedness and bideal properties of compact bilinear operators on asymmetric normed spaces.

In order to discuss closedness we introduce  the extended quasi-metrics $d$ and $d_s$ on  $$\L^2_{p_1,p_2;q}(X,Y;Z)\quad\mbox{and}\quad \L^2_{p_1,p_2;q^s}(X,Y;Z)$$
 given  by
\begin{equation}\label{def1.ext-metric}
d(T_1,T_2)=\sup\{q(T_2(x,y)-T_1(x,y)):(x,y)\in B_{p_1}\times B_{p_2}\}\,,
\end{equation}
and  by
\begin{equation}\label{def2.ext-metric}
d_s(T_1,T_2)=\sup\{q^s(T_2(x,y)-T_1(x,y)):(x,y)\in B_{p_1}\times B_{p_2}\}\,,
\end{equation}
for $T_1,T_2$ in $\L^2_{p_1,p_2;q}(X,Y;Z) $ or in $\L^2_{p_1,p_2;q^s}(X,Y;Z), $  respectively. Notice that $d_s$ is actually  an extended metric  and that
$d$ and $d_s$ are nothing else that the extended quasi-metrics \eqref{def.ext-dist} corresponding to $\|\cdot|_{p_1,p_2;q}$ and $\|\cdot|_{p_1,p_2;q^s}$, respectively.

\begin{prop}\label{p1.closed-id} Let $(X,p_1), (Y,p_2)$ and   $(Z,q)$ be    asymmetric normed spaces. Let $(T_n)$ be a sequence of $q^s$-precompact bilinear operators in $\L^2_{p_1,p_2;q^s}(X,Y;Z)$ and let  $T:X\times Y\to Z$  be a mapping . If for every $(x,y)\in X\times Y$ the sequence  $(T_n(x,y))$ is $q^s$-convergent to $T(x,y)$,  uniformly on $B_{p_1}\times B_{p_2}$,   then $T$ is a bilinear $q^s$-precompact operator.  In particular, the cone $\mathcal K^2_{p_1,p_2;q^s}(X,Y;Z)$ is $d_s$-closed in $\L^2_{p_1,p_2;q^s}(X,Y;Z)$.
\end{prop}\begin{pf}  The bilinearity of $T$ follows from the pointwise convergence of the sequence $(T_n)$, the bilinearity of each $T_n$ and the fact that $q^s$ is a norm on $Z$.

Indeed,
\begin{align*}
  q^s(T(\alpha x,y)-\alpha T(x,y))&\le q^s(T(\alpha x,y)- T_n(\alpha x,y))\\&+q^s(\alpha T_n(x,y)-\alpha T(x,y))\to 0 \;\mbox{ as }\; n\to \infty\,,
\end{align*}
so that $T(\alpha x,y)=\alpha T(x,y)$. The equality  $T(x+x',y)= T(x,y)+T(x',y)$ follows similarly as well as the linearity with respect to the second variable.

 Let now $\varepsilon>0$. The   uniform convergence  on $B_{p_1}\times B_{p_2}$ implies the existence of $n_0\in \mathbb{N}$ such that
 \begin{equation}\label{eq1.closed-id}
 q^s(T_n(x,y)- T(x,y))\le \varepsilon\,,\end{equation}
 for all $n\ge n_0$ and all $(x,y)\in B_{p_1}\times B_{p_2}$.

 The $q^s$-precompactness of $T_{n_0}$  implies the existence of $(x_1,y_1),\dots,(x_m,y_m)\in B_{p_1}\times B_{p_2}$ such that
 for every $(x,y)\in B_{p_1}\times B_{p_2}$ there exists $k\in\{1,\dots,m\}$ with
 \begin{equation}\label{eq2.closed-id}
 q^s(T_{n_0}(x,y)- T_{n_0}(x_k,y_k))\le \varepsilon\,.\end{equation}

For an arbitrary point  $(x,y)  \in B_{p_1}\times B_{p_2}$ choose  $k\in\{1,\dots,m\}$ according to \eqref{eq2.closed-id}. Then, by \eqref{eq1.closed-id},
\begin{align*}
   q^s(T(x,y)- T(x_k,y_k)) &\le q^s(T(x,y)- T_{n_0}(x,y))\\ &+q^s(T_{n_0}(x,y)- T_{n_0}(x_k,y_k))\\&+ q^s(T_{n_0}(x_k,y_k)- T(x_k,y_k))\le 3\varepsilon\,,
\end{align*}
which shows that $T$ is $q^s$-precompact.
\qed\end{pf}

The next proposition is concerned with the ideal properties of precompact bilinear operators.
\begin{prop}\label{p1.bid} Let $(X,p_1), (Y,p_2)$ and  $(Z,q)$ be   asymmetric normed spaces.
\begin{enumerate}\item[\rm 1.] Let $(Z_1,q_1)$ be another asymmetric normed space and  let $R:Z\to Z_1$  be a  linear operator. If $\,T:X\times Y\to Z$ is $q$-precompact and $R$ is $(q,q_1)$-continuous, or $T$ is
$q^s$-precompact  and $R$ is $(q^s,q_1)$-continuous, then $RT$ is $q_1$-precompact or  $q^s_1$-precompact, respectively.
 The same property holds for the outside  and hereditary $q$-precompactness notions when $R$ is $(q,q_1)$-continuous.
\item[\rm 2.] Let $(X',p_1'), (Y',p_2')$ be two other asymmetric normed spaces, let $S_1:X'\to X$,     $S_2:Y'\to Y$ be    continuous linear operators and let  $S:X'\times Y'\to X\times Y$ be given by
$S(x',y')=(S_1x',S_2y')$ for $(x',y')\in X'\times Y'$. If the bilinear operator  $T:X\times Y\to Z$ is  $q$-precompact ($q^s$-precompact), then $T\circ S:X'\times Y'\to Z$ is a $q$-precompact ($q^s$-precompact)  bilinear operator.
\end{enumerate}\end{prop}\begin{pf}
  1. Let $\gamma \in\{q,q^s\}.$ The property follows from the inequality
  $$
  \gamma((R(T(x,y))-R(T(x_k,y_k))\le \|R|\,\gamma(T(x,y)-T(x_k,y_k))\le  \|R|\,\varepsilon\,,$$
  where $\{(x_1,y_1),\dots,(x_m,y_m)\}\subseteq B_{p_1}\times B_{p_2}$ is an $\varepsilon$-net, according to the definitions of precompactness (Definitions \ref{def.prec} and \ref{def.prec-op}).  Here $\|R|=\|R|_{q,q_1}$ if $\gamma =q$ and  $\|R|=\|R|_{q^s,q_1}$ if $\gamma =q^s$ 

  2. Let $\beta_1,\beta_2>0$ be such that
  $$
  p_1(S_1x')\le \beta_1 p_1'(x')\quad\mbox{and}\quad  p_2(S_2y')\le \beta_2 p_2'(y')\,,$$
  for all $(x',y')\in X'\times Y'$. If $\beta:=\max\{\beta_1,\beta_2\},$ then
  $$
  S(B_{p_1'}\times B_{p_2'})\subseteq \beta(B_{p_1}\times B_{p_2})\,,$$
and
    $$
  (T\circ S)(B_{p_1'}\times B_{p_2'})\subseteq \beta\, T(B_{p_1}\times B_{p_2})\,.$$

  This inclusion shows that $(T\circ S)(B_{p_1'}\times B_{p_2'})$ is $q$-precompact ($q^s$-precompact).
\qed\end{pf}

\subsection{Arens' adjoint of a bilinear operator}
Another adjoint of a bilinear operator was proposed by Arens  \cite{arens51a}, \cite{arens51b}.  Let $X,Y,Z$ be normed spaces and $T:X\times Y\to Z$ a continuous bilinear operator.
Define $T^\star:Z^*\times X\to Y^*$ by
\begin{equation}\label{def.Arens-adj}
T^\star(\psi,x)(y)=\psi(T(x,y)),\quad y\in Y\,,
\end{equation}
for $\psi\in Z^*$ and $x\in X$.
\begin{prop}\label{p1.Arens} The so defined adjoint is a continuous bilinear operator with
$$
\|T^\star\|=\|T\|\,.   $$
\end{prop}

The definition can be extended to the asymmetric case  defining $T^\diamond:Z^\flat\times Y\to Y^\flat$ by the   formula \eqref{def.Arens-adj} with $\psi\in Z^\flat$ and $x\in X$.

\begin{remark}
I do not know whether a Schauder type theorem holds or not for this adjoint of a compact bilinear operator, even in the classical case (i.e., for $X,Y,Z$ Banach spaces).
\end{remark}

\textbf{Funding resources.}  This research did not receive any specific grant from funding agencies in the public, commercial, or not-for-profit sectors.

\textbf{Acknowledgements.}  The author thanks the learned referee for the pertinent remarks and suggestions that led to a major improvement of the presentation.

\end{document}